\documentclass[11pt]{article}
\usepackage{amsmath}
\usepackage{amssymb}
\usepackage{vatola}
\def\q{\quad}
\def\qq{\qquad}
\def\t{\text}
\def\qtq#1{\q\t{#1}\q}
\def\f{\frac}
\def\e{\equiv}
\def\a{\alpha}
\def\b{\binom}

\def\m{\pmod}

\def\sls#1#2{(\f{#1}{#2})}
 
\def\Ls#1#2{\Big(\f{#1}{#2}\Big)}

\let \pro=\proclaim
\let \endpro=\endproclaim

\begin{document}
\leftline{Revised: April 19, 2012}
\par\q\par\q
 \centerline {\Large\bf
On the further properties of $\{U_n\}$}
$$\q$$\centerline{Zhi-Hong Sun} $$\q$$
\centerline {School of Mathematical Sciences, Huaiyin Normal
University,} \centerline {Huaian, Jiangsu 223001, P.R. China}
\centerline {E-mail: zhihongsun@yahoo.com}
 \centerline {Homepage:
http://www.hytc.edu.cn/xsjl/szh}
 \abstract{Let $\{U_n\}$ be given by $U_0=1$ and $U_n=-2\sum_{k=1}^{[n/2]}
 \b n{2k}U_{n-2k}\ (n\ge 1)$, where $[\cdot]$ is the greatest
  integer function. In the paper we
present a summation formula and several congruences involving
$\{U_n\}$.
\par\q
}
 \endabstract
 \footnotetext[1] {The author is
supported by the Natural Sciences Foundation of China (grant no.
10971078).}

\section*{1. Introduction}
\par\q\   The Euler numbers $\{E_n\}$ and Euler polynomials $\{E_n(x)\}$
are defined by
$$\align &E_0=1,\q E_n=-\sum_{k=1}^{[n/2]}\binom n{2k}E_{n-2k}\q(n\ge 1),
\\&E_n(x)=\f 1{2^n}\sum_{k=0}^{[n/2]}\b n{2k}(2x-1)^{n-2k}E_{2k},\endalign
 $$ where $[x]$ is the greatest integer not exceeding $x$. In [7]
the author introduced and studied the sequence $\{U_n\}$ (similar to
Euler numbers) as below:
$$U_0=1,\q
U_n=-2\sum_{k=1}^{[n/2]}\b n{2k}U_{n-2k}\q(n\ge 1).\tag 1.1$$ Since
$U_1=0$, by induction we have $U_{2n-1}=0$ for $n\ge 1$.
 The first few values of $U_{2n}$ are shown below:
$$\align &U_2=-2,\q U_4=22,\q U_6=-602,\q U_8=30742,\q U_{10}=-2523002,
\\&U_{12}=303692662,\q U_{14}=-50402079002, \q U_{16}=11030684333782.
\endalign$$
\par Let
$\sls ap$ be the Legendre symbol. In [7], the author proved that for
any prime $p>3$,
$$\sum_{k=1}^{[2p/3]}\f{(-1)^{k-1}}k\e 3p\Ls p3U_{p-3}\pmod{p^2}.\tag 1.2$$
The Bernoulli numbers $\{B_n\}$ and Bernoulli polynomials
$\{B_n(x)\}$ are given by
$$B_0=1,\ \sum_{k=0}^{n-1}\b nkB_k=0\ (n\ge 2)\qtq{and}
B_n(x)=\sum_{k =0}^n\b nkB_kx^{n-k}\ (n\ge 0).$$
  By [7, p.217],
$$B_{p-2}\Ls 13\e 6U_{p-3}\pmod p\qtq{for any prime} p>3.\tag 1.3$$
 In [8] S. Mattarei and R.
Tauraso proved that for any prime $p>3$,
$$\sum_{k=0}^{p-1}\b{2k}k\e \Ls p3-\f{p^2}3B_{p-2}\Ls 13\pmod
{p^3}.$$ Thus,
$$\sum_{k=0}^{p-1}\b{2k}k\e \Ls p3-2p^2U_{p-3}
\pmod{p^3}\qtq{for any prime} p>3.\tag 1.4$$ Suppose that $p$ is a
prime of the form $3k+1$ and so $4p=L^2+27M^2$ with $L,M\in\Bbb Z$
and $L\e 1\pmod 3$. From (1.3) and [3, Theorem 6] we have
$$\b{\f{2(p-1)}3}{\f{p-1}3}
\e \Big(-L+\f pL+\f{p^2}{L^3}\Big)(1+p^2U_{p-3})\e -L+\f
pL+p^2\Big(\f 1{L^3}-LU_{p-3}\Big)\pmod{p^3}.\tag 1.5$$

\par In Section 2 we prove a summation formula involving ${U_n}$,
see Theorem 2.1. Let $\Bbb N$ be the set of positive integers.
 If $n\in\Bbb N$ and
$2^{\a}\mid n$, in [7] the author determined $U_{2n}\pmod
{2^{\a+7}}$.  In Section 3 we prove
$$3U_{2n}\e -3072n^4+4608n^3+2240n^2+1680n+2
\pmod{2^{\a+14}}\qtq{for}n\ge 7.$$
 For $k,m,b\in\Bbb N$ with $2\mid b$, in Section 3 we
also show that
$$U_{2^mk+b}\e U_b+2^{b+1}\pmod{2^{\t{min}\{b,m\}+3}}.$$
Let $k,m\in\Bbb N$ and $b\in\{0,2,4,\ldots\}$. From [7, Theorem 4.3]
we have $U_{k\varphi(3^m)+b}\e U_b\pmod {3^m}$, where $\varphi(n)$
is Euler's totient function. In Section 4 we prove a congruence for
$U_{k\varphi(3^m)+b}-U_b\pmod {3^{m+4}}$ for $m\ge 3$, see Theorem
4.1. In Section 5 we prove a congruence for
$E_{k\varphi(3^m)+b}-(3^b+1)E_b\pmod {3^{m+4}}$ for $m\ge 3$, see
Theorem 5.1.

\section*{2. A summation formula involving $U_n(x)$}
\par For $n=0,1,2,\ldots$ let
$$U_n(x)=\sum_{r=0}^n\b nrU_rx^{n-r}=\sum_{k=0}^{[n/2]}\b
n{2k}U_{2k}x^{n-2k}.\tag 2.1$$ The first few $U_n(x)$ are given
below:
$$\align &U_0(x)=1,\q U_1(x)=x,\q U_2(x)=x^2-2,
\\&U_3(x)=x^3-6x,\q U_4(x)=x^4-12x^2+22,\
\\&U_5(x)=x^5-20x^3+110x,\q U_6(x)=x^6-30x^4+330x^2-602.
\endalign$$
By [7, Theorem 2.3] we have
$$\align &U_n(x-1)-U_n(x)+U_n(x+1)=x^n,\tag 2.2
\\&U_n(x)+U_n(x+3)=(x+1)^n+(x+2)^n,\tag 2.3
\\&U_n(x+3)-U_n(x-3)=(x+2)^n+(x+1)^n-(x-1)^n-(x-2)^n.\tag 2.4\endalign$$
Taking $a_n=U_n(x)$ and $b_n=x^n$ in [7, Theorem 2.2]  we obtain
$$x^n=2\sum_{k=0}^{[n/2]}\b n{2k}U_{n-2k}(x)-U_n(x).$$
That is,
$$U_n(x)=x^n-2\sum_{k=1}^{[n/2]}\b n{2k}U_{n-2k}(x).\tag 2.5$$
Since
$$\align \int_a^bU_n(x)dx&=\sum_{k=0}^n\b nkU_k\int_a^bx^{n-k}dx
=\sum_{k=0}^n\b nkU_k\f{x^{n-k+1}}{n-k+1}\Big|_a^b
\\&=\f 1{n+1}\sum_{k=0}^{n+1}\b{n+1}kU_kx^{n+1-k}\Big|_a^b,
\endalign$$
we see that
$$\int_a^b U_n(x) dx=\f{U_{n+1}(b)-U_{n+1}(a)}{n+1}.\tag 2.6$$
This together with (2.4) yields
$$\int_{a-3}^{a+3}U_n(x) dx=\f{(a+2)^{n+1}+(a+1)^{n+1}
-(a-1)^{n+1}-(a-2)^{n+1}}{n+1}.\tag 2.7$$ Since $U_n(0)=U_n$, by
(2.6) we have
$$U_n(x)=U_n+n\int_0^xU_{n-1}(t)dt.\tag 2.8$$
\par Let $m,n\in\Bbb N$. From [1] we have the following well known summation formulas.
$$\sum_{k=0}^{m-1}k^n=\f{B_{n+1}(m)-B_{n+1}}{n+1}\qtq{and}
\sum_{k=0}^{m-1}(-1)^k k^n=\f{E_n(0)-(-1)^mE_n(m)}2.\tag 2.9$$ Now
we present the following similar result.
 \pro{Theorem 2.1} Let $m,n\in\Bbb N$ and
$$S_n(m)=(m-1)^n+(m-2)^n-(m-4)^n-(m-5)^n+(m-7)^n+(m-8)^n-\cdots,$$
where the term $a^n$ vanishes  when $a\le 0$. Then
$$S_n(m)=\cases U_n(m)-(-1)^{\f m3}U_n&\t{if $3\mid m$,}
\\U_n(m)-(-1)^{[\f{m+1}3]}U_n/2&\t{if $3\nmid m$ and $2\mid n$,}
\\U_n(m)-(-1)^{[\f m3]}U_n(1)&\t{if $3\nmid m$ and $2\nmid n$.}
\endcases$$
\endpro
Proof. Using (2.3) we see that
$$\align
&(m-1)^n+(m-2)^n-(m-4)^n-(m-5)^n+(m-7)^n+(m-8)^n
\\&\q-\cdots-(-1)^{[\f m3]}\Big(\big(m-3[\f m3]+2\big)^n+
\big(m-3[\f m3]+1\big)^n\Big)
\\&=(U_n(m)+U_n(m-3))-(U_n(m-3)+U_n(m-6))+(U_n(m-6)+U_n(m-9))
\\&\q-\cdots -(-1)^{[\f m3]}\Big(U_n\big(m-3[\f
m3]+3\big)+U_n\big(m-3[\f m3]\big)\Big)
\\&=U_n(m)-(-1)^{[\f m3]}U_n\big(m-3[\f m3]\big).
\endalign$$
Thus,
$$S_n(m)=\cases U_n(m)-(-1)^{\f m3}U_n(0)&\t{if $3\mid m$,}
\\U_n(m)-(-1)^{[\f m3]}U_n(1)&\t{if $3\mid m-1$,}
\\(-1)^{[\f m3]}\cdot 1+U_n(m)-(-1)^{[\f m3]}U_n(2)&\t{if $3\mid m-2$.}
\endcases$$
\par Clearly $U_n(0)=U_n$. By (2.3) and (2.1),  we have $U_n(-1)+U_n(2)=1$ and
so $U_n(2)=1-U_n(-1)=1-(-1)^nU_n(1)$. If $2\mid n$, using (1.1) we
see that
$$U_n(1)=\sum_{k=0}^{n/2}\b n{2k}U_{2k}
=\sum_{k=0}^{n/2}\b n{2k}U_{n-2k} =U_n-\f 12U_n=\f 12U_n\tag 2.10$$
and so
$$U_n(2)=1-U_n(1)=1-\f 12U_n.\tag 2.11$$
\par Now putting all the above together we deduce the result.
\par\q
\pro{Corollary 2.1} For $m\in\Bbb N$ we have
$$\aligned
&S_2(m)=\cases m^2-2+2(-1)^{m/3}&\t{if $3\mid m$,}
\\m^2-2+(-1)^{[(m+1)/3]}&\t{if $3\nmid m$,}
\endcases
\\&S_3(m)=\cases m^3-6m&\t{if $3\mid m$,}
\\m^3-6m+5(-1)^{[m/3]}&\t{if $3\nmid m$,}
\endcases
\\&S_4(m)=\cases m^4-12m^2+22(1-(-1)^{m/3})&\t{if $3\mid m$,}
\\m^4-12m^2+11(2-(-1)^{[(m+1)/3]})&\t{if $3\nmid m$.}
\endcases
\endaligned$$
\endpro

\pro{Corollary 2.2}  For $n\in\Bbb N$ we have
$$\align U_{2n}&=\f
23\Big\{2^{2n}+3^{2n}-\sum_{k=1}^n\b{2n}{2k}4^{2k}U_{2n-2k}\Big\}
\\&=\f
23\Big\{7^{2n}+6^{2n}-4^{2n}-3^{2n}+1-\sum_{k=1}^n\b{2n}{2k}
8^{2k}U_{2n-2k}\Big\}.\endalign$$
\endpro
Proof. Taking $m=4,8$ in Theorem 2.1 and replacing $n$ with $2n$ we
see that
$$3^{2n}+2^{2n}=U_{2n}(4)+U_{2n}/2$$
and $$7^{2n}+6^{2n}-4^{2n}-3^{2n}+1=U_{2n}(8)+U_{2n}/2.$$ Since
$$U_{2n}(x)=\sum_{r=0}^n\b{2n}{2r}U_{2r}x^{2n-2r}
=U_{2n}+\sum_{k=1}^n\b{2n}{2k}x^{2k}U_{2n-2k},$$ from the above we
deduce the result.

\section*{3. A congruence for $U_{2n}\pmod {2^{14}}$}
\par Suppose $n\in\{3,4,5,\ldots\}$. From [7, Theorem 4.1 and
Corollary 4.1] we know that
$$U_{2n}\e -16n-42\pmod{2^7}.\tag 3.1$$ Moreover, if $n$ is even and
$2^{\a}\mid n$, then
$$U_{2n}\e 48n+\f 23\pmod{2^{\a+7}}.\tag 3.2$$
\par Let $p$ be a prime and let $\t{ord}_pm$ be the greatest integer
$\a$ such that $p^{\a}\mid m$. If $p^s\le n< p^{s+1}$, then
$$\t{ord}_pn!=\big[\f np\big]+\big[\f n{p^2}\big]+\cdots+ \big[\f
n{p^s}\big]<\f np+\f n{p^2}+\cdots+\f n{p^s}+\cdots =n\cdot \f
{1/p}{1-1/p}=\f n{p-1}.\tag 3.3$$ \pro{Lemma 3.1} Suppose $n\in\Bbb
N$, $n\ge 5$ and $2^{\a}\mid n$. Then
$$\align 3U_{2n}+2^7n(2n-1)U_{2n-2}
&\e 2(7^{2n}+6^{2n}-4^{2n}-3^{2n}+1) +2^{16+2\a}(n-1)\\&\q-23\cdot
2^{13}n(n-1)+7\cdot 2^{15}n(n-1)^3 \pmod{2^{\a+19}}.\endalign$$
\endpro
Proof.  For $3\le k\le n-1$, using (3.3) and (1.1) we have
$8^{2k}/k\e 0\pmod{2^{18}}$ and $2\mid U_{2n-2k}$. For $k=n>3$ we
have $8^{2k}/k\e 0\pmod{2^{19}}$. Thus, for $k\ge 3$,
$$2\b{2n}{2k}8^{2k}U_{2n-2k}=2n\b{2n-1}{2k-1}\f{8^{2k}}kU_{2n-2k}
\e 0\pmod{2^{\a+20}}.$$ Hence, by Corollary 2.2 we get
$$3U_{2n}\e 2\Big\{7^{2n}+6^{2n}-4^{2n}-3^{2n}+1-\b {2n}2
8^2U_{2n-2}-\b{2n}48^4U_{2n-4}\Big\}\pmod{2^{\a+20}}.\tag 3.4$$
Since $U_{2n-4}\e -16(n-2)-42=-16n-10\pmod {2^7}$  and
$$2\b{2n}48^4=2^{12}n(n-1)\f{4(n-1)^2-1}3\e 0\pmod{2^{\a+12}},$$
we see that
$$\align 2\b{2n}48^4U_{2n-4}&
\e 2^{12}n(n-1)\f{4(n-1)^2-1}3(-16n-10)
\\&=-2^{16}n^2(n-1)\f{4(n-1)^2-1}3-2^{13}n(n-1)\f{5(4(n-1)^2-1)}3
\\&\e -2^{16+2\a}(n-1)\cdot 3(4(n-1)-1)-23\cdot
2^{13}n(n-1)(4(n-1)^2-1)
\\&\e -2^{16+2\a}\cdot 3(4(n-1)-(n-1))+23\cdot 2^{13}n(n-1)-23\cdot
2^{15}n(n-1)^3
\\&\e -2^{16+2\a}(n-1)+23\cdot 2^{13}n(n-1)-7\cdot
2^{15}n(n-1)^3 \pmod{2^{\a+19}}.\endalign$$ Hence, by (3.4) and the
fact $2\b{2n}28^2U_{2n-2}=2^7n(2n-1)U_{2n-2}$ we deduce the result.
\par\q
\pro{Theorem 3.1} Let $n\in\Bbb N$ with $n\ge 7$ and $2^{\a}\mid n$.
Then
$$3U_{2n}\e -3072n^4+4608n^3+2240n^2+1680n+2
\pmod{2^{\a+14}}.$$
\endpro
Proof. Since $U_{2n-2}\e -16(n-1)-42\pmod{2^7}$, by Lemma 3.1 we get
$$\align &3U_{2n}+2^7n(2n-1)(-16(n-1)-42)
\\&\e 2(7^{2n}+6^{2n}-4^{2n}-3^{2n}+1) -23\cdot
2^{13}n(n-1)\pmod{2^{\a+14}}.\endalign$$ As $2n\ge \a+13$, we have
$6^{2n}\e 4^{2n}\e 0\pmod{2^{\a+13}}$. We also note that
$2^7n(2n-1)(16(n-1)+42)=2^8(16n^3+18n^2-13n)$. Now, from the above
we deduce
$$ 3U_{2n}\e 2^8(16n^3+18n^2-13n)+2(7^{2n}-3^{2n}+1)
-2^{13}n(n-1)\pmod{2^{\a+14}}.\tag 3.5$$ It is clear that
$$\align 7^{2n}&=(1+48)^n=1+\sum_{k=1}^nn\b{n-1}{k-1}\f{48^k}k
\\&\e 1+48\b n1+48^2\b n2+48^3\b n3\pmod{2^{\a+14}}\endalign$$
and
$$\align 3^{2n}&=(1+8)^n=1+\sum_{k=1}^nn\b{n-1}{k-1}\f{8^k}k
\\&\e 1+8\b n1+8^2\b n2+8^3\b n3+8^4\b n4\pmod{2^{\a+14}}.\endalign$$
Thus,
$$\align 7^{2n}-3^{2n}&\e (48-8)\b n1+(48^2-8^2)\b n2+(48^3-8^3)\b n3-8^4\b n4
\\&\e 40n+1120(n^2-n)-768n(n-1)(n-2)-1536n(n-1(n-2)(n-3)
\\&=-2^9\cdot 3n^4+2^8\cdot 33n^3-2^5\cdot 421n^2+6600n
\\&\e -2^9\cdot 3n^4+2^8n^3+2^5\cdot 91n^2-1592n
\pmod{2^{\a+13}}.\endalign$$ This together with (3.5) yields
$$\align 3U_{2n}&\e 2^8(16n^3+18n^2-13n)+2(-2^9\cdot 3n^4+2^8n^3
+2^5\cdot 91n^2-1592n+1) -2^{13}n(n-1)
\\&=-2^{10}\cdot 3n^4+2^9\cdot 9n^3+2^6\cdot 35n^2+1680n+2
\\&=-3072n^4+4608n^3+2240n^2+1680n+2
\pmod{2^{\a+14}}.\endalign$$ This proves the theorem.

\pro{Lemma 3.2} Let $k,m,b\in\Bbb N$ with $2\mid b$. Then
$$U_{2^mk+b}-U_b\e \f{2^{b+1}}9-\f 23\sum_{r=1}^{\f b2-1}
\b b{2r}2^{2r}(U_{2^mk+b-2r}-U_{b-2r})\pmod{2^{m+3}}.$$
\endpro
Proof. From [7, (4.1)] we have
$$U_{2n}=\f 23\Big(1-\sum_{r=1}^n\b{2n}{2r}2^{2r}U_{2n-2r}\Big).$$
Thus,
$$U_b=\f 23\Big(1-\sum_{r=1}^{b/2}\b b{2r}2^{2r}U_{b-2r}\Big)$$ and
 $$U_{2^mk+b}= \f 23\Big\{1-2^{2^mk+b}U_0
-\sum_{r=1}^{2^{m-1}k+\f b2-1} (2^mk+b)\cdots (2^mk+b-2r+1)\cdot
\f{2^{2r}}{(2r)!}U_{2^mk+b-2r}\Big\}.$$  By (3.3), $2^{2r}/(2r)!\e
0\pmod 2$. By (1.1), $2\mid U_{2n}$ for $n\ge 1$. We also have
$2^mk+b\ge m+2$. Thus, from the above we deduce
$$\align U_{2^mk+b}&\e \f 23\Big\{1-\sum_{r=1}^{2^{m-1}k+\f b2-1}
b(b-1)\cdots(b-2r+1)\f{2^{2r}}{(2r)!}U_{2^mk+b-2r}\Big\}
\\&=\f 23\Big\{1-\sum_{r=1}^{b/2}\b b{2r}2^{2r}U_{2^mk+b-2r}\Big\}
\pmod{2^{m+3}}.\endalign$$ Therefore,
$$U_{2^mk+b}-U_b\e -\f 23\sum_{r=1}^{b/2}\b b{2r}2^{2r}
(U_{2^mk+b-2r}-U_{b-2r})\pmod{2^{m+3}}.$$ By (3.2), $U_{2^mk}\e
48\cdot 2^{m-1}k+\f 23\e \f 23\pmod{2^{m+3}}$. So we have
$$\align U_{2^mk+b}-U_b&\e -\f{2^{b+1}}3(U_{2^mk}-1)-
\f 23\sum_{r=1}^{\f b2-1}\b b{2r}2^{2r} (U_{2^mk+b-2r}-U_{b-2r})
\\&\e \f{2^{b+1}}9-
\f 23\sum_{r=1}^{\f b2-1}\b b{2r}2^{2r} (U_{2^mk+b-2r}-U_{b-2r})
\pmod {2^{m+3}}.\endalign$$ This is the result.

 \pro{Theorem
3.2} Let $k,m\in\Bbb N$.
\par $(\t{\rm i})$ If $b\in\{2,4,6,\ldots\}$, then
$$U_{2^mk+b}\e U_b+2^{b+1}\pmod{2^{\t{min}\{b,m\}+3}}.$$
\par $(\t{\rm ii})$ We have
$$U_{2^mk+2}\e -\f{10}9\pmod{2^{m+3}}\qtq{and}
U_{2^mk+4}\e \f{34}3\pmod{2^{m+3}}.$$
\par $(\t{\rm iii})$ If $b\in\{4,6,8,\ldots\}$ and $b\le m-2$, then
$$U_{2^mk+b}\e U_b+2^{b+1}(4b+5)\pmod{2^{b+5}}.$$
\endpro
Proof. If $b\in\{2,4,6,\ldots\}$, by Lemma 3.2 we have
$$\aligned &U_{2^mk+b}-U_b-\f{2^{b+1}}9
\\&\e -\f 23\sum_{r=1}^{\f b2-1}\b b{2r}2^{2r}
\big(U_{2^mk+b-2r}-U_{b-2r}-\f{2^{b-2r+1}}9\big) -\f 23\cdot
\f{2^{b+1}}9\sum_{r=1}^{\f b2-1}\b b{2r}
\\&= -\f 23\sum_{r=1}^{\f b2-1}\b b{2r}2^{2r}
\big(U_{2^mk+b-2r}-U_{b-2r}-\f{2^{b-2r+1}}9\big) -\f
{2^{b+2}}{27}(2^{b-1}-2) \pmod {2^{m+3}}.\endaligned\tag 3.6$$
Hence,
$$\align U_{2^mk+b}-U_b-\f{2^{b+1}}9\e
-\f 23\sum_{r=1}^{\f b2-1}\b b{2r}2^{2r}
\big(U_{2^mk+b-2r}-U_{b-2r}-\f{2^{b-2r+1}}9\big) \pmod
{2^{\t{min}\{b,m\}+3}}.\endalign$$ Therefore, for $b=2$,
$U_{2^mk+b}-U_b-\f{2^{b+1}}9\e 0\pmod{2^{\t{min}\{b,m\}+3}}$. Now we
prove (i) by induction on $b$. Suppose that the congruence
$$U_{2^mk+b-2r}-U_{b-2r}-\f{2^{b-2r+1}}9\e 0\pmod{2^{\t{min}\{m,b-2r\}+3}}$$
holds for $r=1,2,\ldots,\f b2-1$. As
$$2\cdot 2^{2r}\cdot 2^{\t{min}\{m,b-2r\}+3}
\e\cases 2\cdot 2^{2r}\cdot 2^{b-2r+3}\e 0\pmod{2^{b+3}}&\t{if $b\le
m$,}
\\2\cdot 2^{m+1}\cdot 2^3\e 0\pmod{2^{m+3}}&\t{if $b\ge m$ and $r>\f
m2$,}
\\2\cdot 2^{2r}\cdot 2^{m-2r+3}\e 0\pmod{2^{m+3}}&\t{if $b\ge m\ge
2r$,}\endcases$$
 from the above and induction we
deduce $U_{2^mk+b}-U_b-\f{2^{b+1}}9\e 0\pmod{2^{\t{min}\{b,m\}+3}}$.
This yields (i).
\par Now we consider (ii). Putting $b=2$ in Lemma 3.2 we see that
$$U_{2^mk+2}\e U_2+\f{2^3}9=-2+\f 89=-\f{10}9\pmod{2^{m+3}}.$$
Taking $b=4$ in Lemma 3.2 and then applying the above we deduce
$$U_{2^mk+4}-U_4\e \f{2^5}9-\f 23\b 42\cdot 2^2(U_{2^mk+2}-U_2)
\e \f{32}9-16\cdot \f 89=-\f{32}3\pmod{2^{m+3}}$$ and so
$U_{2^mk+4}\e U_4-\f{32}3=22-\f{32}3=\f{34}3\pmod{2^{m+3}}$. This
proves (ii).
\par Finally we consider (iii). Assume $2\le b\le m-2$.
 By (i), for  $1\le r\le \f b2-1$ we have
$U_{2^mk+b-2r}-U_{b-2r}-\f{2^{b-2r+1}}9\e 0\pmod {2^{b-2r+3}}$.
Thus, it follows from (3.6) that $$U_{2^mk+b}-U_b-\f{2^{b+1}}9 \e
-\f{2^{b+2}}{27}(2^{b-1}-2)\e 2^{2b+1}-2^{b+3}\m{2^{b+4}}.$$ Using
this and (3.6) we get
$$\aligned &U_{2^mk+b}-U_b-\f{2^{b+1}}9
\\&\e
-\f 23\sum_{r=1}^{\f b2-1}\b
b{2r}2^{2r}(2^{2(b-2r)+1}-2^{b-2r+3})
-\f{2^{b+3}}{27}(2^{b-2}-1)
\\&=-\f{2^{b+2}}3\sum_{r=1}^{\f b2-1}\b b{2r}2^{b-2r}+\f{2^{b+4}}3
\sum_{r=1}^{\f b2-1}\b b{2r}-\f{2^{2b+1}}{27}+\f{2^{b+3}}{27}
\\&\e \cases-\f{2^{b+2}}3\b b{b-2}\cdot 2^2+\f{2^{b+4}}3(2^{b-1}-2)
-2^{b+3}\e 2^{b+3}(b-1)\m{2^{b+5}}&\t{if $b>2$,}
\\0 \m{2^{b+5}}&\t{if $b=2$}\endcases
\endaligned$$
and therefore for $b>2$,
$$U_{2^mk+b}-U_b\e \f{2^{b+1}}9+2^{b+3}(b-1)
\e 2^{b+1}(4b+5) \m{2^{b+5}}.$$ This proves (iii). The proof is now
complete. \pro{Corollary 3.1} Let $k,m,b\in\Bbb N$ with $2\mid b$.
Then $U_{2^mk+b}\e U_b\m{2^{\t{min}\{b,m\}+1}}$.
\endpro
Proof. This is immediate from Theorem 3.2(i).

 \section*{4. A congruence for $U_{k\varphi(3^m)+b}\pmod {3^{m+4}}$}

 \par In [5] the author proved that for $k,m\in\Bbb N$, $m\ge 4$ and
 $b\in\{0,2,4,\ldots\}$,
$$E_{2^mk+b}-E_b\e \cases 5\cdot 2^mk\pmod{2^{m+4}}&\t{if $b\e
0,6\pmod 8$,}\\-3\cdot 2^mk\pmod{2^{m+4}}&\t{if $b\e 2,4\pmod 8$.}
\endcases$$ A generalization to Euler polynomials was given in [6,
Theorem 3.3].

\par From the proof of [7, Theorem 4.2] we have the following
lemma.
 \pro{Lemma 4.1} For $n\in\Bbb N$ we have
 $$2^{2n}U_{2n}=\sum_{k=0}^n\b{2n}{2k}3^{2k}E_{2k}.$$
\endpro

 \pro{Theorem 4.1}
Let $k,m\in\Bbb N$, $m\ge 3$ and $b\in\{0,2,4,\ldots\}$. Then
$$U_{k\varphi(3^m)+b}-U_b
\e\cases 3^mk(9b-40)\pmod{3^{m+4}}&\t{if $3\mid b$,}
\\-3^mk\cdot 22\pmod{3^{m+4}}&\t{if $3\mid b-1$,}
\\-3^mk(9b-32)\pmod{3^{m+4}}&\t{if $3\mid b-2$.}
\endcases$$
\endpro
Proof. By (3.3), $\t{ord}_3(2r)!\le r-1$. Thus, for $r\ge 3$ we have
$2r-\t{ord}_3(2r)!\ge 2r-(r-1)=r+1\ge 4$ and so $3^{2r}/(2r)!\e
0\m{3^4}$.  Clearly $3^2\mid b(b-1)\cdots (b-2r+1)$ for $r\ge 3$.
Thus, for $r\ge 3$,
$$\aligned &(k\varphi(3^m)+b)(k\varphi(3^m)+b-1)\cdots
(k\varphi(3^m)+b-2r+1)
\\& \e b(b-1)\cdots (b-2r+1)+k\varphi(3^m)\sum_{i=0}^{2r-1}
\f{b(b-1)\cdots (b-2r+1)}{b-i}
\\&\e b(b-1)\cdots (b-2r+1)\m{3^m}.\endaligned$$
Hence,
$$\b{k\varphi(3^m)+b}{2r}3^{2r}\e \b b{2r}3^{2r}\m{3^{m+4}}\qtq{for}
r\ge 3.\tag 4.1$$ Since $E_0=1,\ E_2=-1$ and $E_4=5$, using (4.1)
and Lemma 4.1 we see that
$$\align &2^{k\varphi(3^m)+b}U_{k\varphi(3^m)+b}\\&=1-9\b{k\varphi(3^m)+b}2
+5\cdot 3^4\b{k\varphi(3^m)+b}4 +\sum_{r=3}^{(k\varphi(3^m)+b)/2}
\b{k\varphi(3^m)+b}{2r}3^{2r}E_{2r}
\\&\e -\f 92(4\cdot 3^{2m-2}k^2+2\cdot 3^{m-1}k(2b-1))
+\f{3^3\cdot
5}8(k^4\varphi(3^m)^4+(4b-6)k^3\varphi(3^m)^3\\&\q+(6b^2-18b+11)
k^2\varphi(3^m)^2+(4b^3-18b^2+22b-6)k\varphi(3^m))
+\sum_{r=0}^{b/2}\b b{2r}3^{2r}E_{2r}
\\&\e -2\cdot 3^{2m}k^2-(2b-1)3^{m+1}k-(4b^3+4b-6)3^{m+2}k+2^bU_b
\pmod{3^{m+4}}.
\endalign$$
By (3.3), $\t{ord}_3r\le \t{ord}_3r!<\f r2$. Thus, for $r\ge 4$ we
have
$$2r-\t{ord}_3r>2r-\f r2=\f{3r}2>6\qtq{and
so}\varphi(3^{m-1})\f{9^r}r=2\cdot 3^{m-2}\cdot\f{3^{2r}}r\e
0\m{3^{m+4}}.$$ Hence,
$$\align &2^{k\varphi(3^m)}-1\\&
=(1-9)^{k\varphi(3^{m-1})}-1 =\sum_{r=1}^{k\varphi(3^{m-1})}
\b{k\varphi(3^{m-1})}r(-9)^r
\\&=\b{k\varphi(3^{m-1})}1(-9)+\b{k\varphi(3^{m-1})}2(-9)^2
+\b{k\varphi(3^{m-1})}3(-9)^3\\&\q+\sum_{r=4}^{k\varphi(3^{m-1})}
k\varphi(3^{m-1})\b{k\varphi(3^{m-1})-1}{r-1}\f{(-9)^r}r
\\&\e
-9k\varphi(3^{m-1})+81\f{k\varphi(3^{m-1})(k\varphi(3^{m-1})-1)}2
+\f{k\varphi(3^{m-1})(k\varphi(3^{m-1})-1)(k\varphi(3^{m-1})-2)}6(-9)^3
\\&\e 3^mk(16+2\cdot 3^mk)\pmod{3^{m+4}}.\endalign$$
Thus,
$$\align &2^{k\varphi(3^m)+b}U_{k\varphi(3^m)+b}-2^bU_b
\\&\e (1+3^mk(16+2\cdot 3^mk))2^bU_{k\varphi(3^m)+b}-2^bU_b
\\&=2^b(U_{k\varphi(3^m)+b}-U_b)+3^mk(16+2\cdot 3^mk)
2^bU_{k\varphi(3^m)+b}\pmod{3^{m+4}}.\endalign$$ By Lemma 4.1,
$$2^{2n}U_{2n}\e E_0+\b{2n}23^2E_2=1-9n(2n-1)\pmod{81}.$$ Thus,
$$2^{k\varphi(3^m)+b}U_{k\varphi(3^m)+b}\e 1-\f
92(k\varphi(3^m)+b)(k\varphi(3^m)+b-1)\e 1-9\b b2\pmod{81}$$ and so
$$\align 2^bU_{k\varphi(3^m)+b}&\e\f{1-9\b b2}{2^{k\varphi(3^m)}}
\e \f{1-\f 92b(b-1)}{1+3^mk(16+2\cdot 3^mk)}
\\&\e\big(1-\f 92b(b-1)\big)(1-3^mk(16+2\cdot 3^mk))
\\&\e 1-\f 92b(b-1)-16\cdot 3^mk
\pmod{81}.\endalign$$ Therefore,
$$\aligned &-2\cdot 3^{2m}k^2-(2b-1)3^{m+1}k-(4b^3+4b-6)3^{m+2}k
\\&\e 2^{k\varphi(3^m)+b}U_{k\varphi(3^m)+b}-2^bU_b
\e 2^b(U_{k\varphi(3^m)+b}-U_b)+3^mk(16+2\cdot 3^mk)
2^bU_{k\varphi(3^m)+b}
\\&\e 2^b(U_{k\varphi(3^m)+b}-U_b)+3^mk(16+2\cdot 3^mk)
\big(1-\f 92b(b-1)-16\cdot 3^mk\big)
\\&\e\cases
 2^b(U_{k\varphi(3^m)+b}-U_b)+3^mk(16+9b(b-1))\pmod{3^{m+4}}&\t{if $m\ge 4$,}
\\2^b(U_{k\varphi(3^3)+b}-U_b)+27k(16+9b(b-1))+3^6k^2
\pmod{3^7}&\t{if $m=3$.}\endcases\endaligned$$ This yields
$$ 2^b(U_{k\varphi(3^m)+b}-U_b)
\e -3^mk(36b^3+9b^2+33b-41)\pmod{3^{m+4}}.$$ If $3\mid b$, then
$2^{-b}=(1-9)^{-\f b3}$ and so
$$\align U_{k\varphi(3^m)+b}-U_b
&\e (1-9)^{-\f b3}(-3^mk)(36b^3+9b^2+33b-41)\\& \e
-(1+3b)(33b-41)3^mk \e (9b-40)3^mk\pmod{3^{m+4}}.\endalign$$ If
$3\mid b-1$, then $b^2\e 2b-1\pmod 9$,
 $b^3\e 1\pmod 9$ and $2^{-b}=4\cdot 2^{-b-2}=4(1-9)^{-\f{b+2}3}$. Thus,
$$\align U_{k\varphi(3^m)+b}-U_b
&\e 4(1-9)^{-\f {b+2}3}(-3^mk)(36b^3+9b^2+33b-41)\\& \e
-4(1+3(b+2))(36+9(2b-1)+33b-41)3^mk
\\&=-4(3b+7)(51b-14)3^mk\e 8(3b+7)(15b+7)3^mk\\&\e 8(45(2b-1)+126b+49)3^mk
=8(216b+4)3^mk\\&\e 8(216+4)3^mk\e -22\cdot
3^mk\pmod{3^{m+4}}.\endalign$$ If $3\mid b-2$, then $b^2\e
-2b-1\pmod 9$, $b^3\e -1\pmod 9$ and $2^{-b}=2\cdot
2^{-b-1}=-2(1-9)^{-\f{b+1}3}$. Thus,
$$\align U_{k\varphi(3^m)+b}-U_b
&\e -2(1-9)^{-\f {b+1}3}(-3^mk)(36b^3+9b^2+33b-41)\\& \e
2(1+3(b+1))(-36+9(-2b-1)+33b-41)3^mk
\\&=10(3b+4)(3b-1)3^mk\e 10(9(-2b-1)+9b-4)3^mk\\&\e 10(-9b-13)3^mk
\e -(9b-32)3^mk\pmod{3^{m+4}}.\endalign$$ This completes the proof.

\section*{5. A congruence for $E_{k\varphi(3^m)+b}\pmod {3^{m+4}}$}
 \pro{Lemma 5.1} For $n\in\Bbb N$ we have
$$(3^{2n}+1)E_{2n}=\sum_{r=0}^n\b{2n}{2r}2^{2n-2r+1}3^{2r}E_{2r}.$$
\endpro
Proof. By [6, Theorem 2.1 and Lemma 2.1],
$$\f
12(3^{2n}+1)E_{2n}=\sum_{r=0}^n\b{2n}{2r}(1-3)^{2n-2r}3^{2r}E_{2r}.$$
This is the result.
\par Let $k,m\in\Bbb N$ and $b\in\{0,2,4,\ldots\}$. From [2, p. 231]
or [4, Corollary 7.1] we have
$$E_{k\varphi(3^m)+b}\e (3^b+1)E_b\m{3^m}.$$
\par Now we prove the following stronger congruence.
 \pro{Theorem 5.1} Let $k,m\in\Bbb N$, $m\ge 3$
and $b\in\{0,2,4,\ldots\}$. Then
$$E_{k\varphi(3^m)+b}-(3^b+1)E_b\e \cases (9b+20)3^mk\m{3^{m+4}}
&\t{if $3\mid b$,}\\-16\cdot 3^mk \m{3^{m+4}} &\t{if $3\mid
b-1$,}\\(-9b+11)3^mk\m{3^{m+4}} &\t{if $3\mid b-2$.}\endcases$$
\endpro
Proof. As $\varphi(3^m)\ge m+4$, using Lemma 5.1 and (4.1) we see
that
$$\align &E_{k\varphi(3^m)+b}\\&\e
(3^{k\varphi(3^m)+b}+1)E_{k\varphi(3^m)+b}
=\sum_{r=0}^{(k\varphi(3^m)+b)/2}\b{k\varphi(3^m)+b}{2r}
2^{k\varphi(3^m)+b-2r+1}3^{2r}E_{2r}
\\&=2^{k\varphi(3^m)+b+1}E_0+\b{k\varphi(3^m)+b}2
2^{k\varphi(3^m)+b-1}3^2E_2+\b{k\varphi(3^m)+b}4
2^{k\varphi(3^m)+b-3}3^4E_4\\&\qq+\sum_{r=3}^{b/2}\b
b{2r}2^{k\varphi(3^m)+b-2r+1}3^{2r}E_{2r}\m{3^{m+4}}.
\endalign$$
From the proof of Theorem 4.1 we know that $3^{2r}/(2r)!\e 0\m{3^4}$
for $r\ge 3$. By Euler's theorem, $2^{k\varphi(3^m)}\e 1\m{3^m}$.
Thus, from the above we deduce
$$\align &E_{k\varphi(3^m)+b}\\&\e
2^{k\varphi(3^m)+b+1}
-9(k\varphi(3^m)+b)(k\varphi(3^m)+b-1)2^{k\varphi(3^m)+b-2}
+\b{k\varphi(3^m)+b}42^{b-3}\cdot 81\cdot 5
\\&\qq+\sum_{r=3}^{b/2}\b
b{2r}2^{b-2r+1}3^{2r}E_{2r}
\\&\e 2^{b+1}(2^{k\varphi(3^m)}-1)-9(k^2\varphi(3^m)^2+(2b-1)k\varphi(3^m)+b(b-1))
2^{k\varphi(3^m)+b-2}+9b(b-1)2^{b-2} \\&\qq
+\Big\{\b{k\varphi(3^m)+b}4-\b b4\Big\}2^{b-3}\cdot 81\cdot 5
+\sum_{r=0}^{b/2}\b b{2r}2^{b-2r+1}3^{2r}E_{2r}
\\&\e
2^{b+1}(2^{k\varphi(3^m)}-1)-2^bk^2\cdot
3^{2m}-9(2b-1)2^{b-2}k\varphi(3^m)-9b(b-1)2^{b-2}(2^{k\varphi
(3^m)}-1)\\&\qq+5\cdot 3^3\cdot
2^{b-6}(4b^3-18b^2+22b-6)k\varphi(3^m)+(3^b+1)E_b
\\&\e 2^{b-2}(8-9b(b-1))(2^{k\varphi(3^m)}-1)-2^bk^2\cdot
3^{2m}-(2b-1)2^{b-1}k\cdot 3^{m+1}\\&\qq-2^{b-2}(2b^3+2b-3)k\cdot
3^{m+2}+(3^b+1)E_b
 \m{3^{m+4}}.
\endalign$$
By the proof of Theorem 4.1,
$$2^{k\varphi(3^m)}-1\e 16k\cdot 3^m+2k^2\cdot 3^{2m}\m{3^{m+4}}.$$
Thus,
$$\align E_{k\varphi(3^m)+b}-(3^b+1)E_b&\e
2^{b-2}3^mk\{16(8-9b^2+9b)-6(2b-1)-9(2b^3+2b-3)\}
\\&\qq+2^{b-2}(8-9b(b-1))\cdot 2k^23^{2m}-2^bk^23^{2m}
\\&\e 2^{b-2}(-18b^3+18b^2-48b-1)3^mk\m{3^{m+4}}.
\endalign$$
If $3\mid b$, then $2^{b-2}\e -20(1-9)^{\f b3}\e -20(1-3b)\m{3^4}$.
Thus,
$$\align &E_{k\varphi(3^m)+b}-(3^b+1)E_b\\&\e
-20(1-3b)(-18b^3+18b^2-48b-1)3^mk \e 20(1-3b)(1+48b)3^mk  \\&\e
20(1+45b)\e (9b+20)3^mk \m{3^{m+4}}.\endalign$$ If $3\mid b-1$, then
$b^2\e 2b-1\m 9$, $b^3\e 1\m 9$ and
$$2^{b-2}=-\f 12(1-9)^{\f{b-1}3}\e -\f
12\Big(1-9\cdot\f{b-1}3\Big)=\f 12(3b-4)\m{3^4}.$$ Thus,
$$\align &E_{k\varphi(3^m)+b}-(3^b+1)E_b\\&\e
\f 12(3b-4)(-18b^3+18b^2-48b-1)3^mk \e \f
12(3b-4)(-18+18(2b-1)-48b-1)3^mk
\\&\e (3b-4)(-6b+22)3^mk\e (-18b^2+9b-7)3^mk
\e (-18(2b-1)+9b-7)3^mk\\&=(-27b+11)3^mk\e -16\cdot
3^mk\m{3^{m+4}}.\endalign$$ If $3\mid b-2$, then $b^2\e -2b-1\m 9$,
$b^3\e -1\m 9$ and $2^{b-2}=(1-9)^{\f{b-2}3}\e 1-9\cdot
\f{b-2}3=7-3b\m {3^4}$. Thus,
$$\align &E_{k\varphi(3^m)+b}-(3^b+1)E_b\\&\e
(7-3b)(-18b^3+18b^2-48b-1)3^mk \e (7-3b)(18+18(-2b-1)-48b-1)3^mk
\\&\e (3b-7)(3b+1)3^mk\e (9(-2b-1)-18b-7)3^mk
\e (-9b+11)3^mk\m{3^{m+4}}.\endalign$$ This completes the proof.


\begin{thebibliography}{99}
\bibitem{1} H. Bateman, Higher transcendental functions, Vol.I,
McGraw-Hill Book Co. Inc., 1953.
\bibitem{2} K. W. Chen, Congruences for Euler numbers, Fibonacci
Quart. 42(2004), 128-140.
\bibitem{3} J. B. Cosgrave and K. Dilcher, Mod $p^3$ analogues of
theorems of Gauss and Jacobi on binomial coefficients, Acta Arith.
142(2010), 103-118.
\bibitem{4} Z. H. Sun, Congruences involving Bernoulli polynomials,
Discrete Math. 308(2008), 71-112.
\bibitem{5} Z. H. Sun, Euler numbers modulo $2^n$, Bull. Austral.
Math. Soc. 82(2010), 221-231.
\bibitem{6} Z. H. Sun, Congruences for sequences similar to Euler
numbers, J. Number Theory 132(2012), 675-700.
\bibitem{7} Z. H. Sun, Identities and congruences for a new
sequence, Int. J. Number Theory 8(2012), 207-225.
\bibitem{8} S. Mattarei and R. Tauraso, Congruences for central
binomial sums and finite polylogarithms, preprint, arXiv:1012.1308.
http://arxiv.org/abs/1012.1308.

\end{thebibliography}
\end{document}